# Investigation of bifurcations in the behavior of the process equation

Fahimeh Nazarimehr [a], Sajad Jafari [a,*], Seyed Mohammad Reza Hashemi Golpayegani [a], Louis H. Kauffman[b]

[a] Biomedical Engineering Department, Amirkabir University of Technology, Tehran 15875-4413, Iran

[b] Department of Mathematics, Statistics and Computer Science, University of Illinois, Chicago, IL 60607-7045, USA

**Abstract:** This paper investigates the different behaviors of the process equation and parameters of their occurrences. The process equation is a multistable one dimensional map with nonlinear feedback and can show various behaviors such as period doubling route to chaos, bios, unstable windows and periodic windows. In this note, we focus on different behaviors of the process equation by a deep look at phase portraits and cobweb plots. Therefore, period doubling route to chaos and unifurcations of the equation are investigated, and also the parameter of its entrance to biotic pattern is discussed. The control parameter for the process equation is g (the coupling constant). In higher $g$, the system shows some unstable and periodic windows among the biotic behaviors. Different patterns of these windows and the reasons of their happenings are investigated mathematically. In addition, some other types of unstable and periodic windows are discussed and Q-curves determine the difference of these windows with the previous ones.

## 1. Introduction

The process equation is a one dimensional map with a sine nonlinearity and a single parameter $g$ [1-3]. This equation can generate different behaviors such as convergence to $\pi$, a cascade of bifurcations, chaos, bios and infinitation with respect to changing $g$ parameter. In other words, the recursion goes to the fixed point $A = \pi$ (or other odd multiples of $\pi$ depending on initial conditions) in small $g$ near zero. Increasing $g$ leads to a period doubling route to chaos that is very similar to the bifurcation diagram of logistic map until parameter $g = 4.605$, near the Feigenbaum constant [1-3]. Then, biotic phase starts which has many similarities to biological time series, especially to heartbeat data [1, 2, 4-8]. Bios is an expansive pattern and has been proposed as a mechanism for the expansion of the universe [9]. In addition, biotic patterns can model many structures [10] such as Schrödinger's wave function [11], Klein- Gordon- Fock equation [12], air

* Corresponding author.

*E-mail address:* sajadjafari@aut.ac.ir (S. Jafari).

and ocean temperature and Nile floods [13], sequences of bases in DNA [13], population dynamics [14], prices and stocks [15], and musical compositions [16]. The process equation is a model for creative processes [17] and embodies two basic concepts of process theory [1]. The linear flow information is modeled in iterations and the revolution of complementary opposites is modeled by positive and negative feedback of the sine function [1].

Many researches discuss about various phases of process equation in the viewpoint of its time series, cobweb plot (graph of equation), recurrence plot, and complement plot [1, 2, 10, 18-26]. another map that is derived from the process equation is a circle map that is $mod\ P$ of the process equation and cannot show biotic pattern due to the mixing that occurs from the modulus. Another derived map from the process equation is the Sabelli attractor which is discussed in [2].

In this paper, we aim to analyze different behaviors of the process equation with a deep look to its phase space. In the following section, we describe different phases of the process equation's behavior such as period doubling route to chaos, bios, unstable windows and periodic windows. Discussion about the system behaviors is considered in section 3. And finally, the conclusion is presented in section 4.

## 2. Behaviors of process equation

Bios is a rare type of organization containing the ability of creation and physiological regulation [2, 27]. It is proposed as a mathematical pattern that is a variation of chaos and identified in heartbeat intervals [1, 6, 18]. The process equation, which is shown in Eq. 1, is a map with rare bifurcation that can show evolution of life from one fixed point to chaos, bios and infinity [2].

$$A_{t+1} = f(A_t) = A_t + g\sin(A_t) \quad (1)$$

The process equation is a simple iterative equation that uses a sine function as a harmonic feedback to model fundamental properties of nature [1]. Different behaviors of the equation with respect to the parameter $g$ are described in previous works [1, 2, 18]. In this section, we aim to investigate different behaviors of process equation by a deep look to the phase space variations of the iterative map and its cobweb plots.

### 2.1. Period doubling route to chaos

The behavior of process equation varies when the feedback gain $g$ increases. Fig. 1 shows the bifurcation of this system with respect to the $g$ parameter. As it is shown a period doubling route to chaos occurs by increasing $g$ until $g = 4.6$. In other words, this part of the bifurcation diagram has the same behaviors as the bifurcation diagram of the logistic map except for the appearance of unifurcations. By increasing the parameter $g$ from 0 to 4.6, gain of the sine feedback increases. The system has an infinite number of equilibria located on the multiples of $\pi$. In the interval $g \in$ [0,2], fixed points located at the odd multiples of $\pi$ are stable because the amplitude of the sine is not too high and we have,

| | |
|---|---|
| $f'(A_t) = 1 + g\cos(A_t) \rightarrow \begin{cases} f'(2k\pi) = 1 + g \\ f'((2k+1)\pi) = 1 - g \end{cases}.$ | (2) |

Thus, fixed points located at odd multiples of $\pi$ are stable in the interval $|1-g| < 1$, but fixed points located at even multiples of $\pi$ are unstable in $|1+g| > 1$. Increasing $g$ over 2 makes the odd multiples of $\pi$ unstable and process equation does not have any stable fixed points. Simultaneously, a doubled number of fixed points in the second iterate of the process map ($x_{n+2}$ versus $x_n$) tend to be stable just as the logistic map, and a period doubling route to chaos happens until $g = 4.6$. Fig. 2 shows the cobweb plots of the equation with three different initial conditions in the phase space by changing the gain of the feedback $g$. Fig. 2 investigates that increasing feedback gain $g$ causes increasing in the amplitude of sine and varies the state of the system to be more complex. In addition, Eq. 1 has infinitely many stable intervals which can be trapped onto each of them depended upon initial conditions. In other words, this system is a multistable one, with an infinite number of stabilities. Fig. 2 shows three of these stabilities in three different colored cobwebs. As the $g$ parameter reaches $\pi$, the process equation has a wonderful behavior which is called unifurcation [18]. This behavior appeares missing one limb of its period doubling route to chaos form. The reason of unifurcation is a kind of multistability in the higher iterate of process map. Consider part (a) of Fig. 3 which shows fixed points of the second iterate of the process map in the interval $A \in [0,2\pi)$. Stable fixed points are shown in black color and unstable ones are shown in red. As the figure shows, the two stable fixed points in the second iterate of process equation goes to instability in $g = \pi$ and another four stable fixed points are created. Thus, there are two multistable periodic two behavior in the interval $g \in [\pi, 3.445]$ and the state of the map can switch between these two behaviors depending on its initial condition. They can create two higher limb or two lower ones as shown in bifurcation diagrams of Fig. 4. Part (b) of Fig. 3 shows the fourth iterate of process equation in the interval $A \in [0,2\pi)$. The figure shows that the process equation has a multistable behavior in the periodic four behavior since it has eight stable fixed points in the interval $g \in [3.445,3,5]$. It is obvious that the multistabe behavior follows in the route to chaos (Fig. 1). In the parameter $g$ greater than 4.6, the chaotic behavior changes and a specific pattern occurs which is called Bios. In fact, the amplitude of the sine riding on the identity function increases until the range of one cycle of the sine crosses the range of another sine cycle (Fig. 5) and bifurcation to bios happens. Eq. 3 gives the exact number for the gain of the bifurcation to biotic patterns.

| | |
|---|---|
| $f'(A) = 1 + g\cos(A) \rightarrow f' = 0 \rightarrow A^* = \cos^{-1}(-\frac{1}{g})$<br>$f(A^*) = f(2\pi) \rightarrow A^* + g\ \sin(A^*) = 2\pi \rightarrow g = 4.60333884875170035255658202910\,3$ | (3) |

Fig. 1 part (b) uses the log function to show the biotic and chaotic behavior in the same range. It shows that the amplitude of time series suddenly changes in the bifurcation to biotic behavior.

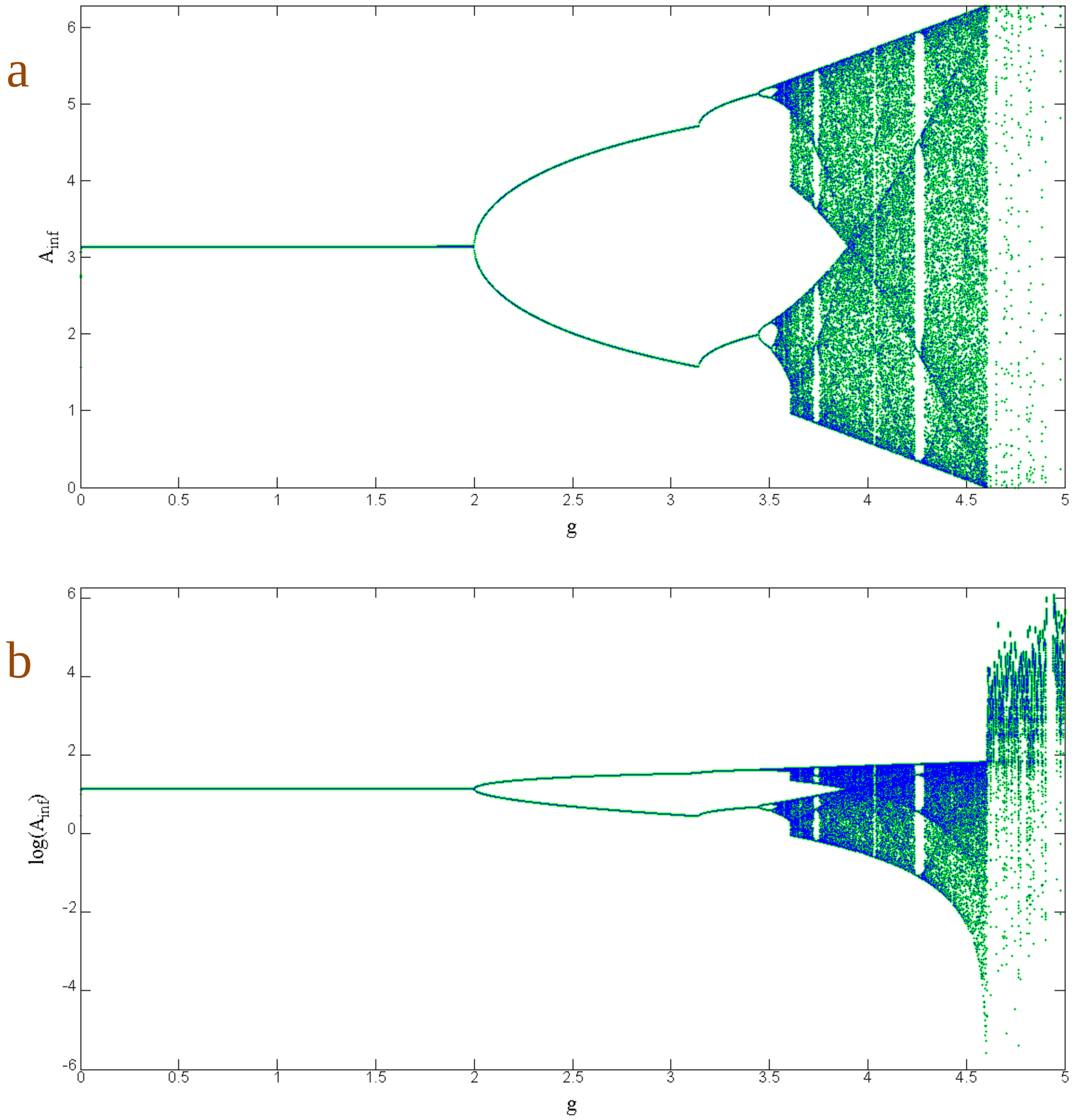

**Fig. 1.** (a) Bifurcation diagram of process equation in the interval $[0,2\pi]$ for $A_{inf}$ with respect to the $g$ parameter and initial condition $A_0 = \frac{\pi}{2}$. (b) Logarithm of the bifurcation diagram of process equation with respect to the $g$ parameter and initial condition $A_0 = \frac{\pi}{2}$.

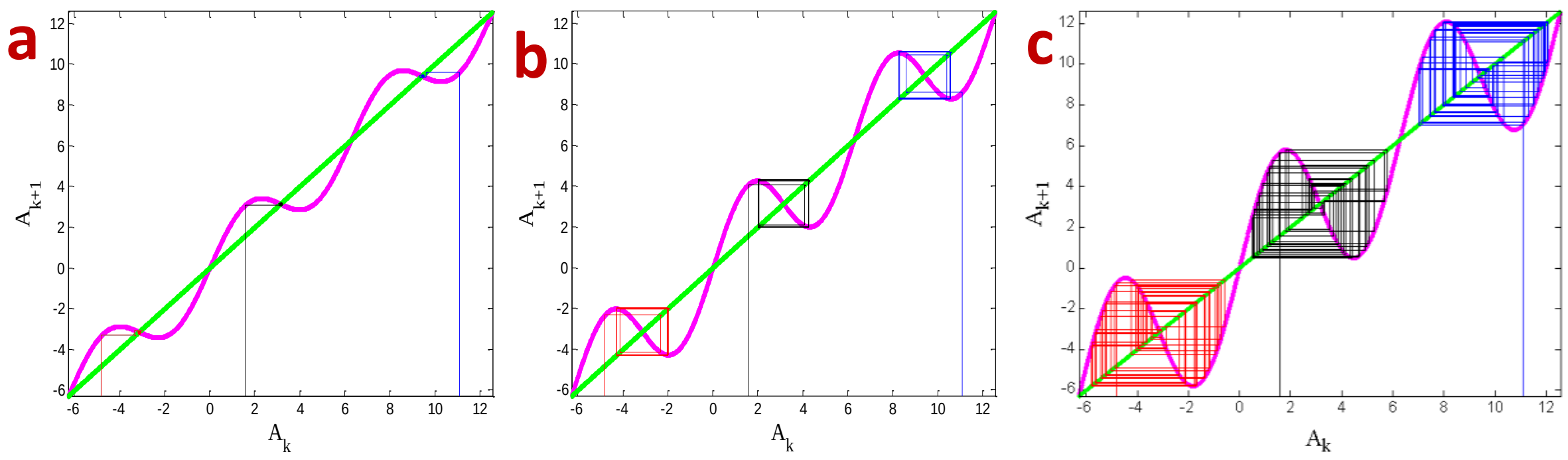

**Fig. 2.** Cobweb plot of process equation in three initial conditions $A_0 = -\frac{3\pi}{2} - 0.1$ (red color), $A_0 = \frac{\pi}{2}$ (black color), and $A_0 = \frac{7\pi}{2} + 0.1$ (blue color). (a) Cobweb of Eq.1 in parameter $g = 1.5$. (b) Cobweb plot of Eq.1 in parameter $g = 2.5$. (c) Cobweb plot of Eq.1 in parameter $g = 4.1$.

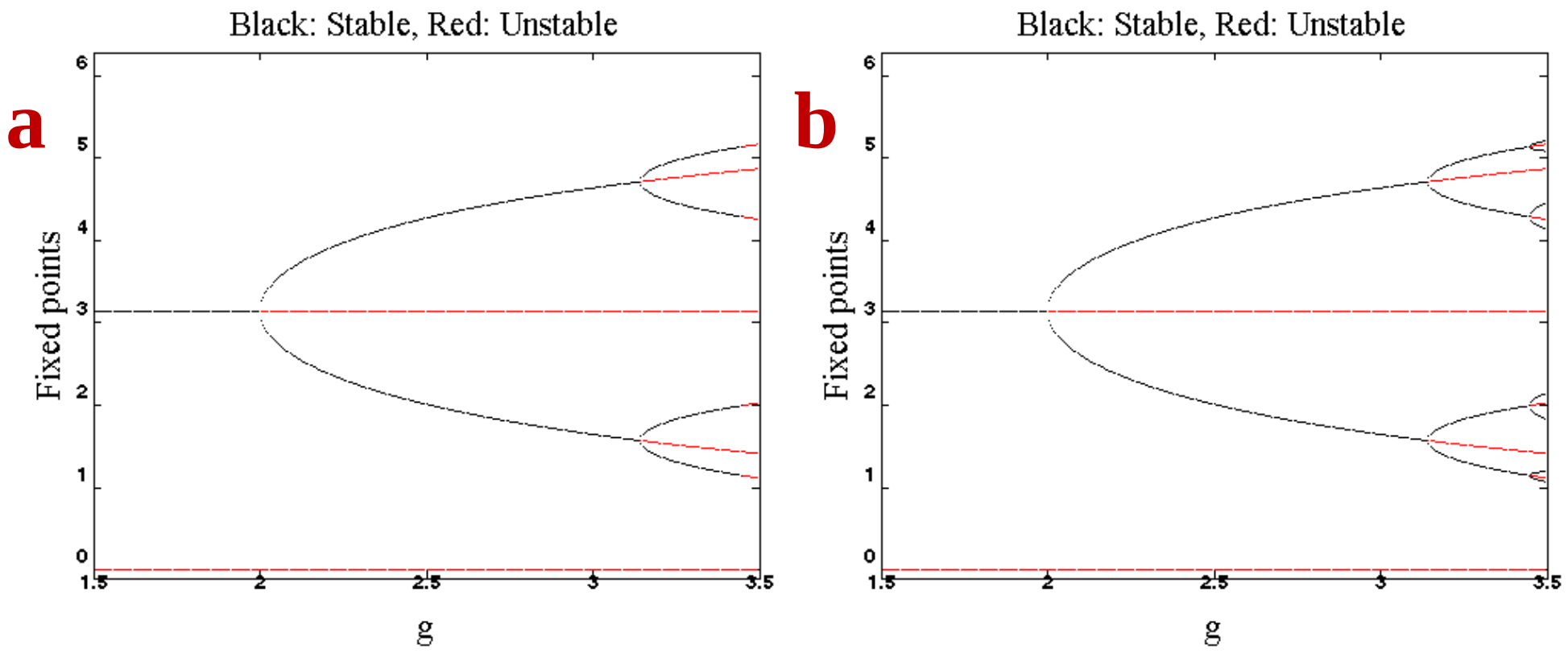

**Fig. 3.** (a) Stable (black) and unstable (red) fixed points of second iterate of process equation in the interval $A \in [0,2\pi)$ with respect to changing $g$ parameter. (b) Stable (black) and unstable (red) fixed points of fourth iterate of process equation in the interval $A \in [0,2\pi)$ with respect to changing $g$ parameter in the interval $g \in [1.5,3.5]$.

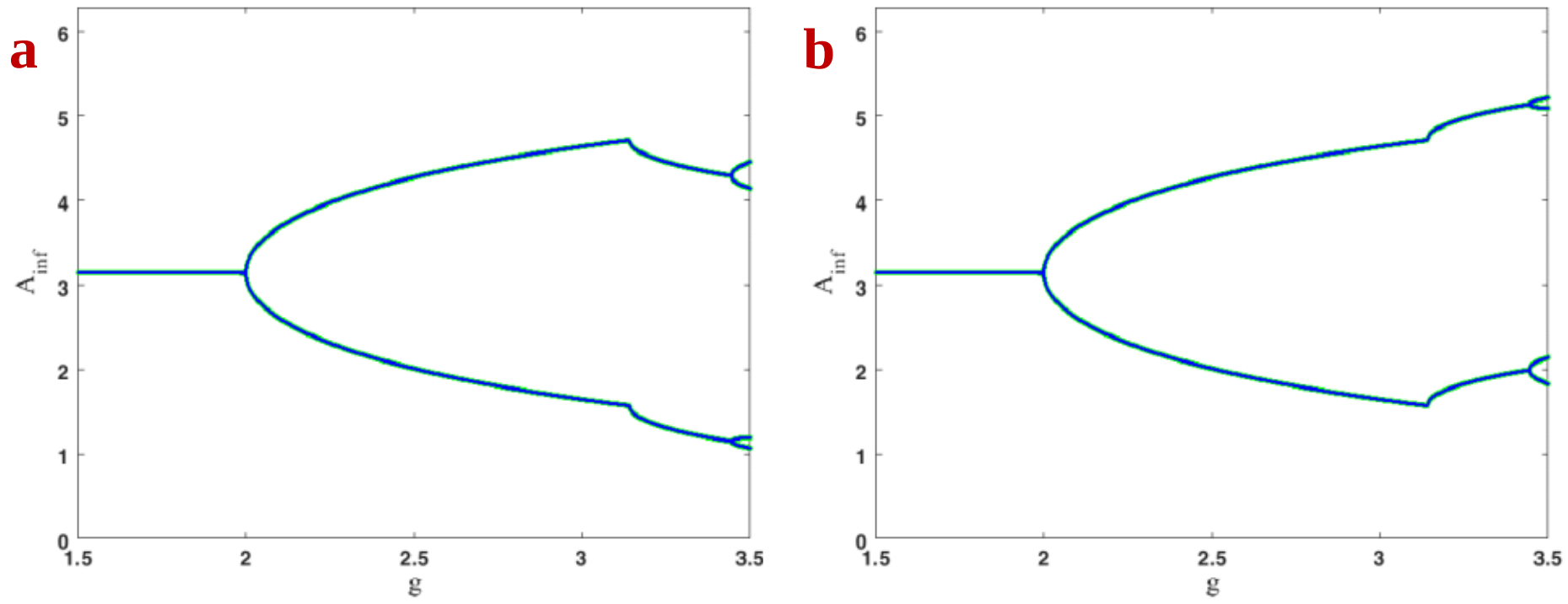

**Fig. 4.** Bifurcation diagram of process equation in the interval $[0,2\pi]$ for $A_{inf}$ with respect to $g$ parameter (a) with initial condition $A_0 = \frac{\pi}{2} - 0.2$. (b) with initial condition $A_0 = \frac{\pi}{2}$.

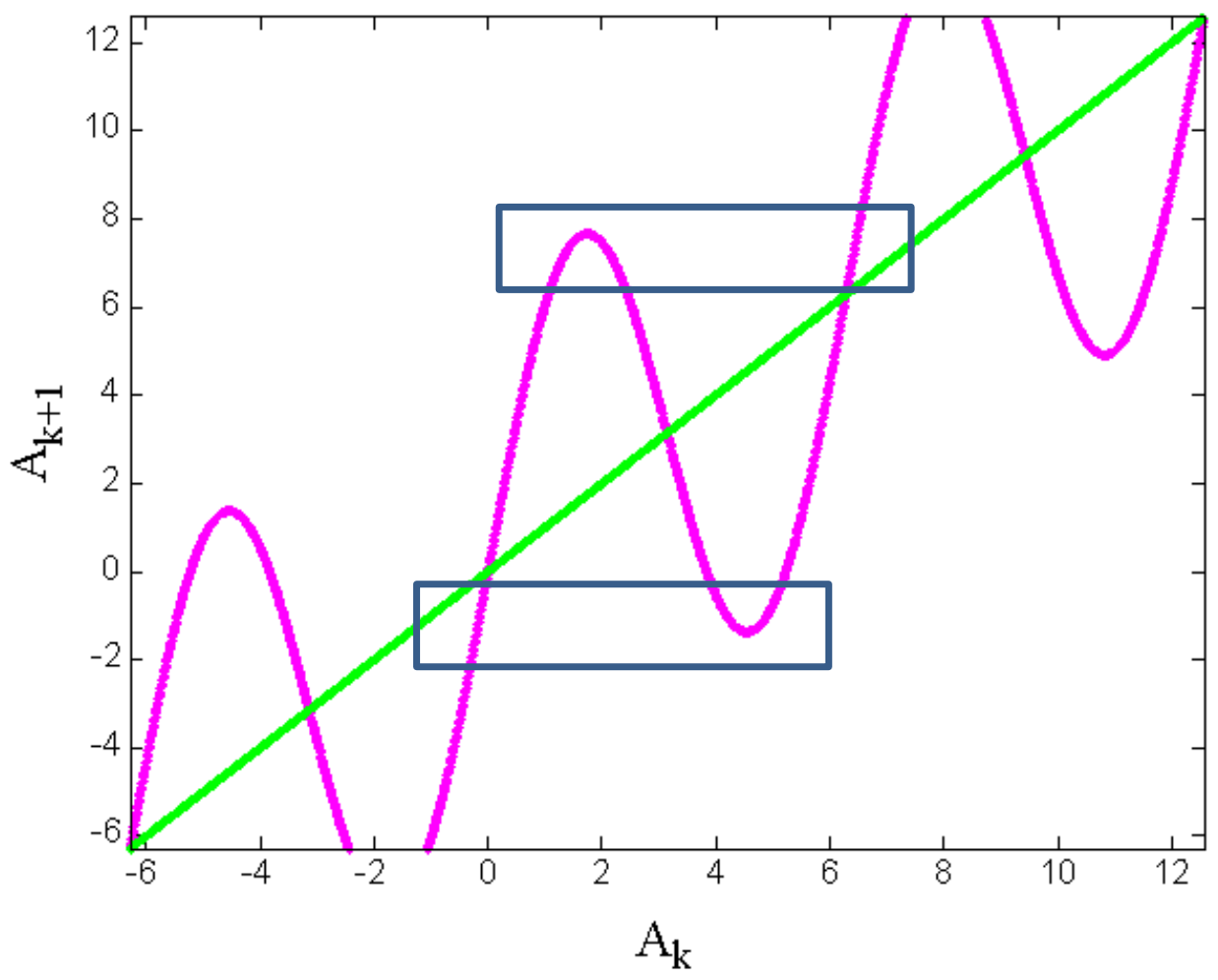

**Fig. 5.** Escaping regions from one cycle of sine to the other one in rectangular intervals.

In the next subsection we investigate properties of biotic patterns.

## 2.2. Biotic behavior

Biotic patterns have wider ranges than chaotic ones because the trajectory can walk through different periods of sine by being in the rectangular intervals of Fig. 5. Fig. 6 shows cobweb plot and time series of biotic behavior in $g = 4.62$. As the figure shows, the amplitude of system state increases by a chaotic walk through different periods of sine function riding on the identity line.

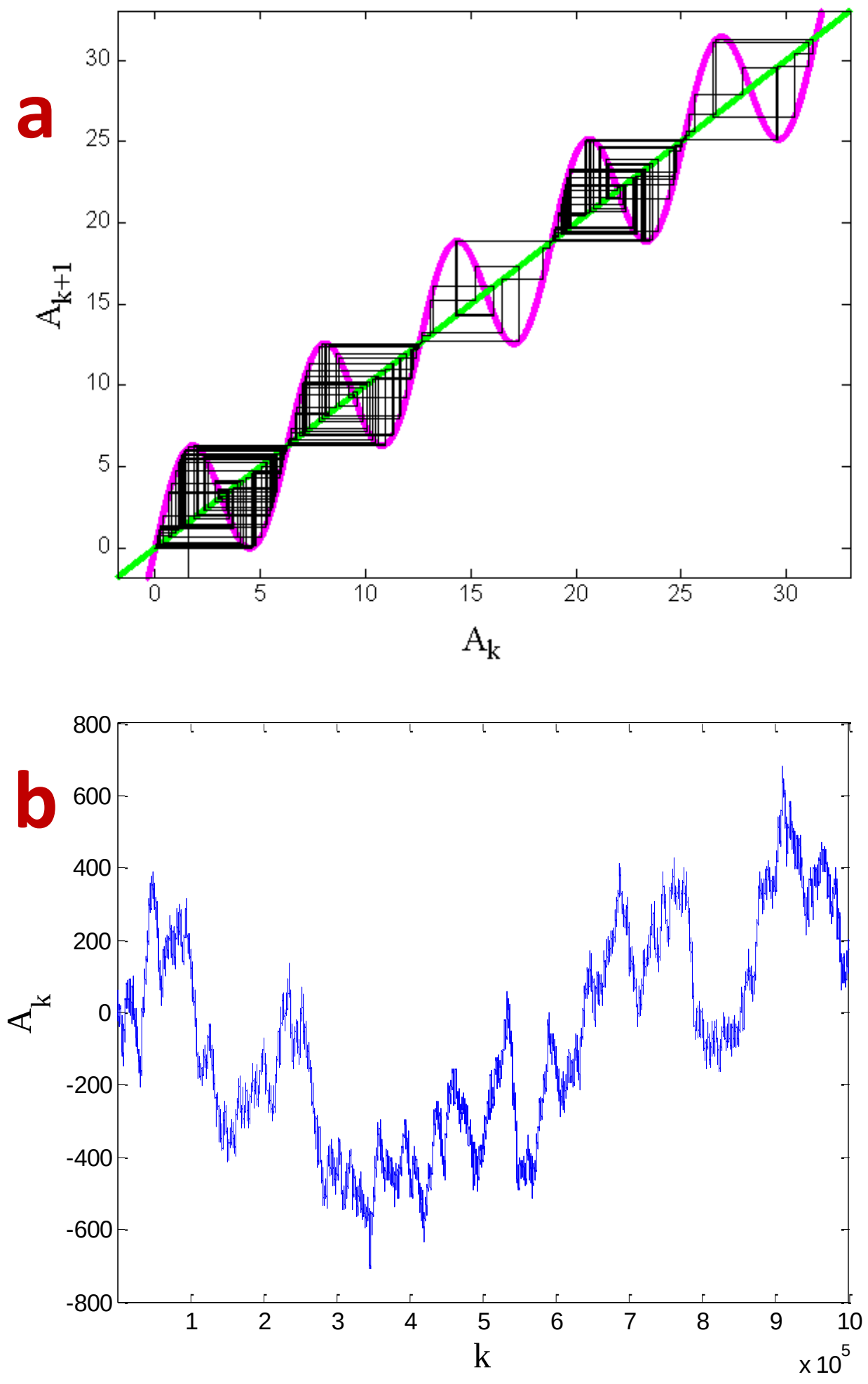

**Fig. 6.** Biotic behavior in parameter $g = 4.62$ and initial condition $A_0 = \frac{\pi}{2}$. (a) Cobweb plot. (b) Time series.

### 2.3. Unstable windows

Another bifurcation happens when $g$ increases and the time series go to infinity. Fig. 7 shows the logarithm (part a) and $mod\ 2\pi$ (part b) of the bifurcation diagram of Eq. 1. Some unstable windows are shown in the black ellipsoids of Fig. 7.

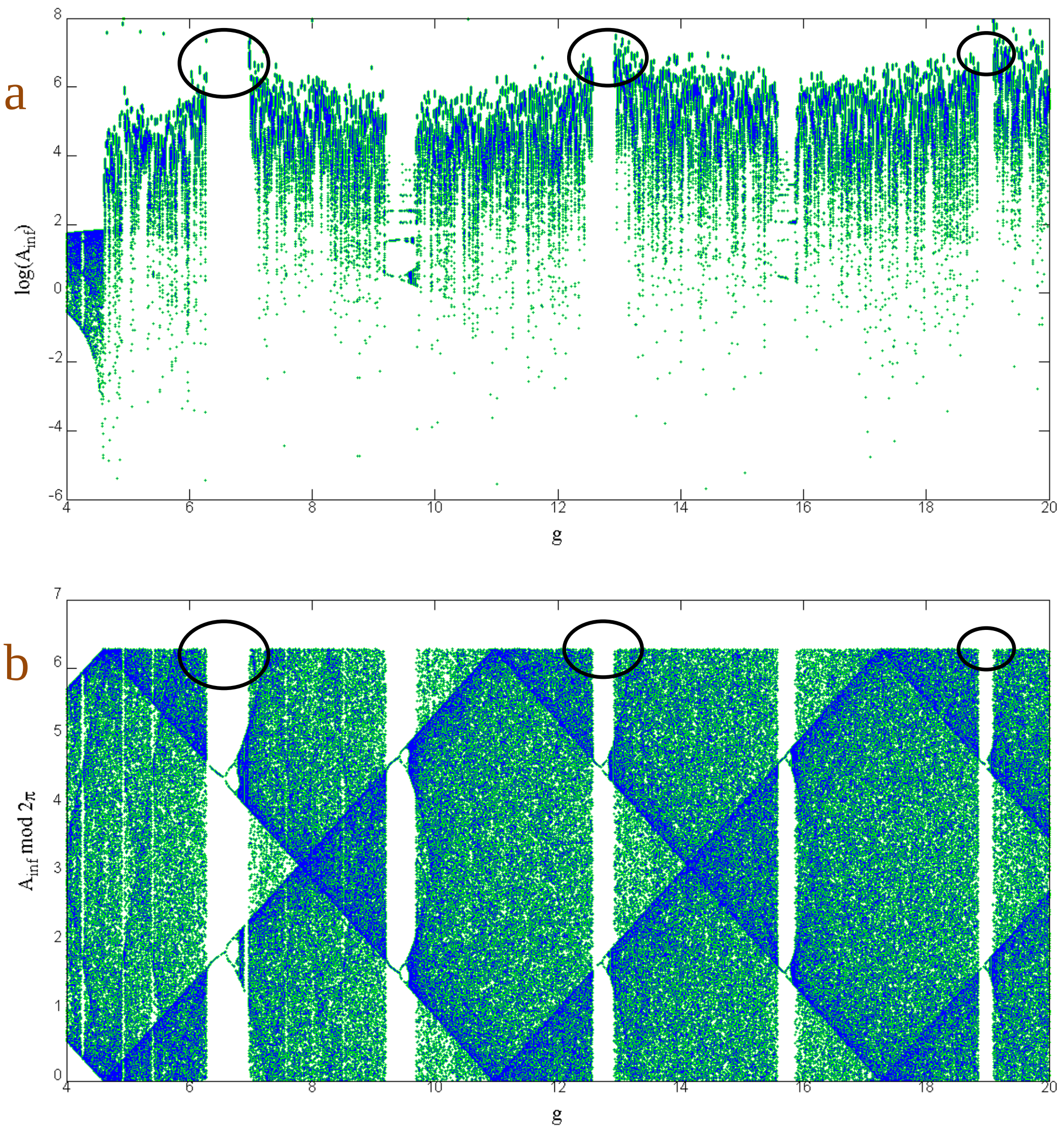

**Fig. 7.** Bifurcation diagram of Eq. 1 with random initial conditions in the interval $A_0 \in [0,2\pi]$. (a) Infinity windows in the bifurcation diagram are shown with ellipsoid signs. (b) Infinity windows in $mod\ 2\pi$ of the bifurcation diagram are shown with ellipsoid black signs.

Fig. 8 shows the phase portrait and cobweb plot of infinity trajectory in $g = 2\pi$. According to this figure, the state of the process equation goes to infinity if the state of $(k+1)$th iteration be in the same phase of sine as the state of $k$th iteration. Thus, the beginning of this type of infinity windows can be calculated using the following condition,

$$\begin{aligned} &A_{k+1} - A_k = 2k\pi \rightarrow \\ &A_k + g * \sin(A_k) - A_k = 2k\pi \rightarrow \\ &g * sin(A_k) = 2k\pi. \end{aligned} \tag{4}$$

Furthermore, the state of system can go to infinity if it has the bios condition and also satisfy the Eq. 4 condition. The minimum $g$ that can satisfy these conditions is the start of first infinity window of this type. Eq. 5 depicts the minimum conditions for $g$ that satisfy the previous conditions and are the beginning of infinity windows.

| | |
|---|---|
| $A_k = \frac{\pi}{2} \rightarrow g = 2k\pi$ | (5) |

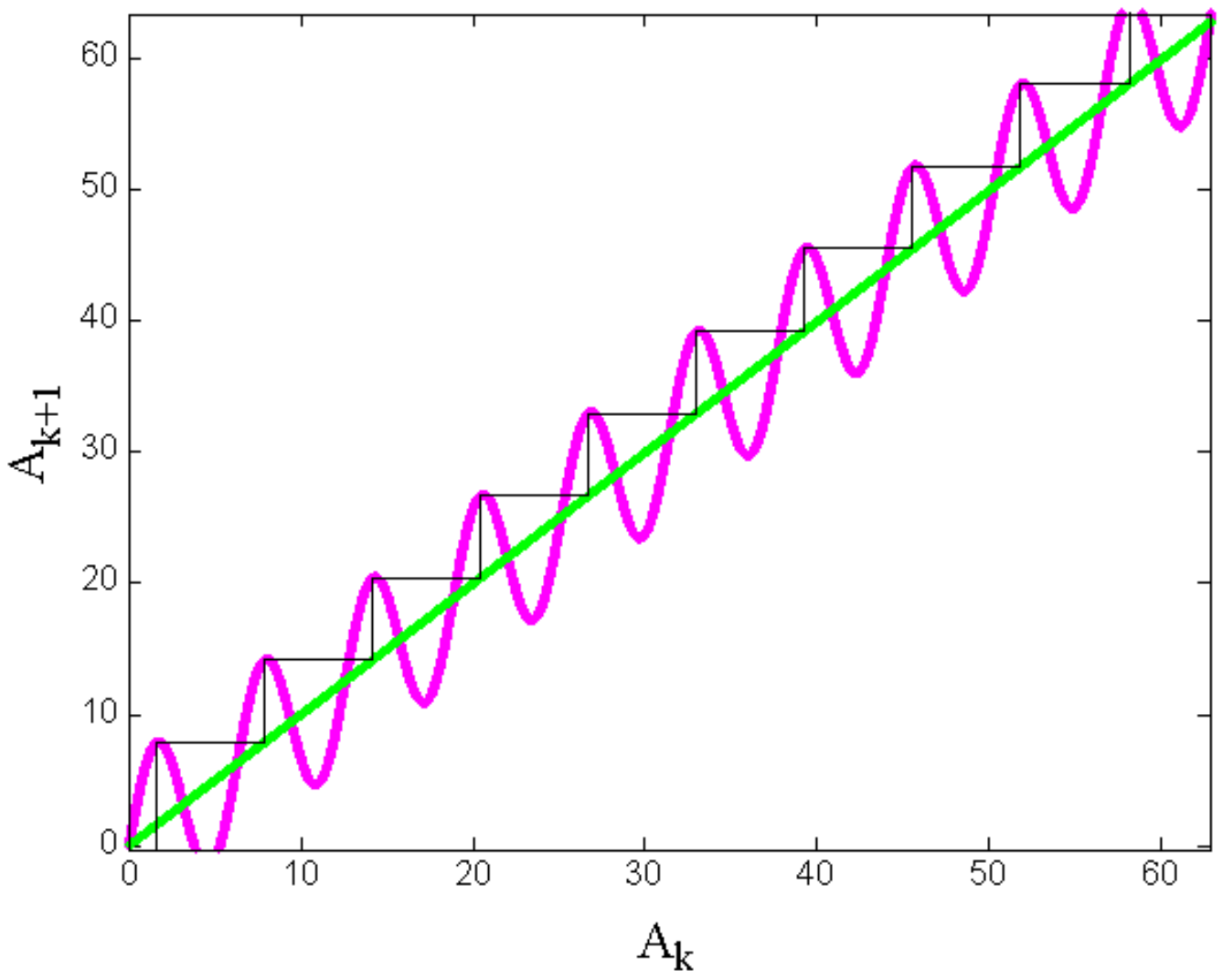

**Fig. 8.** Cobweb plot of system's state in $g = 2\pi$ and initial condition $A_o = \frac{\pi}{2}$.

Now we focus on the first infinity window of this type and try to investigate the system's different behaviors. In the infinity window, there are three categories of dynamical behaviors that we called them periodic fixed point, periodic period $n$ and periodic chaos.

The first category, periodic fixed point, is the condition that state of the system in each iteration has the same phase in the sine portrait as previous iteration. Fig. 8 shows an example of this situation.

The second category, periodic period $n$, is a situation that state of the system in each iteration has the same phase in the cycle of sine with $n$th previous iteration. As an example, consider the periodic period 2 case. Parts (a) and (b) of Fig. 9 show the cobweb plot of the system in periodic period 2 and its phase portrait in $mod\ 2\pi$ (the period of sine function). In this case the state comes back to its phase in the sine after 2 iterations.

The third category, periodic chaos, is a case that the system state never comes back to its phase location in the sine phase and goes to infinity in a chaotic manner. Parts (c) and (d) of Fig. 9 show an example of this situation in $g = 6.8832$. The figure investigates that $mod\ 2\pi$ of the state has a manner just liked the logistic map in going to positive infinity.

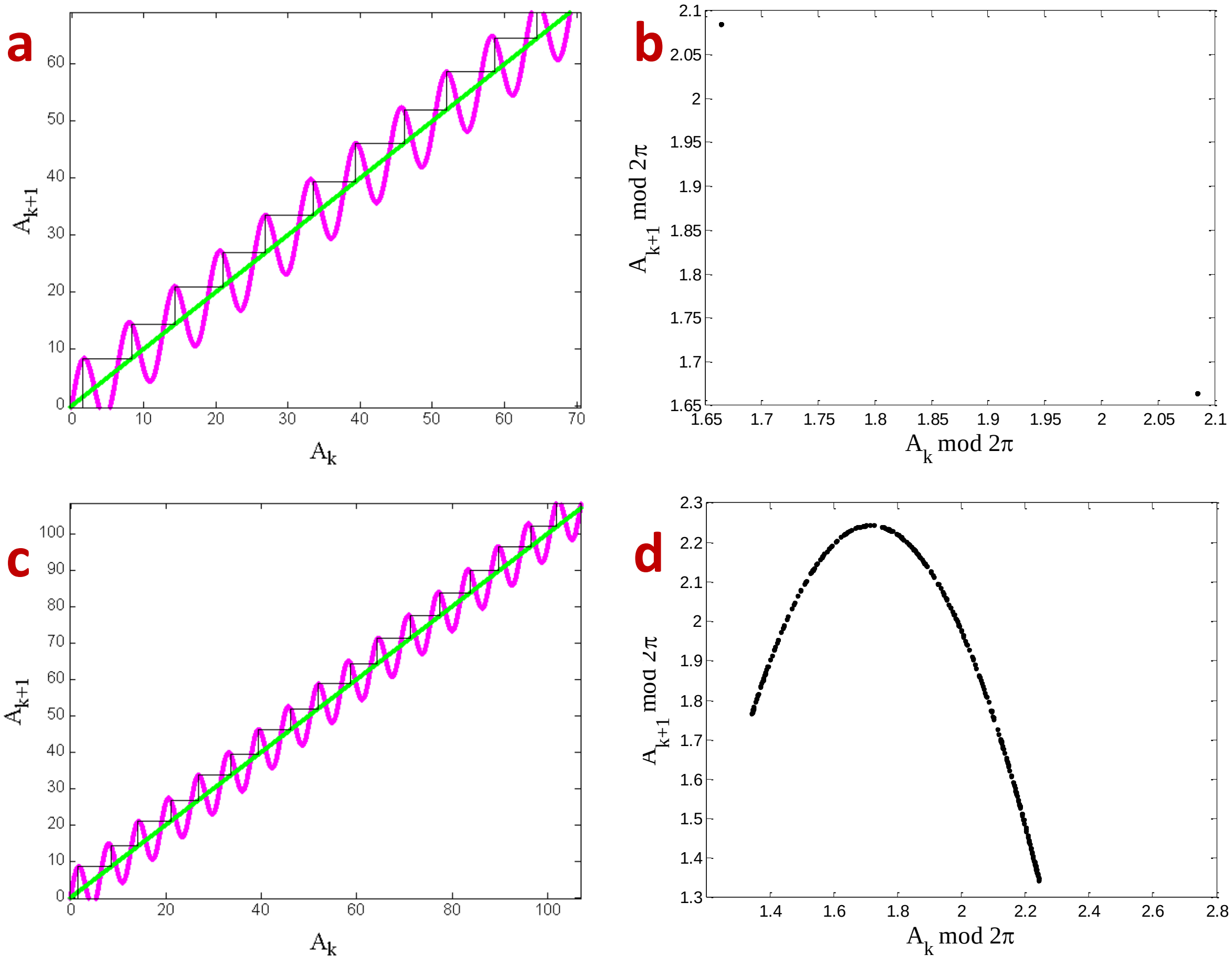

**Fig. 9.** Cobweb plot and $mod\ 2\pi$ phase portrait of Eq.1. (a) Cobweb plot in $g = 6.7332$ and initial condition $A_o = \frac{\pi}{2}$. (b) $Mod\ 2\pi$ phase portrait in $g = 6.7332$ and initial condition $A_o = \frac{\pi}{2}$. (c) Cobweb plot in $g = 6.8832$ and initial condition $A_o = \frac{\pi}{2}$. (d) ) $Mod\ 2\pi$ phase portrait in $g = 6.8832$ and initial condition $A_o = \frac{\pi}{2}$.

The behavior of system in going to infinity has a great sensitivity to initial condition. Part (a) of Fig. 10 proposes bifurcation diagram of process equation in $mod\ 2\pi$ with respect to increasing $g$ parameter in the range of first mentioned unstable window and random initial condition in the interval $A_0 \in [0,2\pi]$. The figure shows that there are two attractors coexisting with each other and a period doubling route to chaos happens to each of them. The lower attractor in the interval $[0,\pi]$ happens when the state goes to positive infinity through the positive half cycle of sine and the upper one in the interval $[\pi, 2\pi]$ happens when the state goes to negative infinity through the negative half cycle of sine. Parts (b) and (c) of Fig. 10 investigate this situation with two constant initial conditions $A_o = \frac{\pi}{2}$ in the positive half cycle of sine and $A_0 = \frac{3\pi}{2}$ in the negative half cycle of sine, respectively.

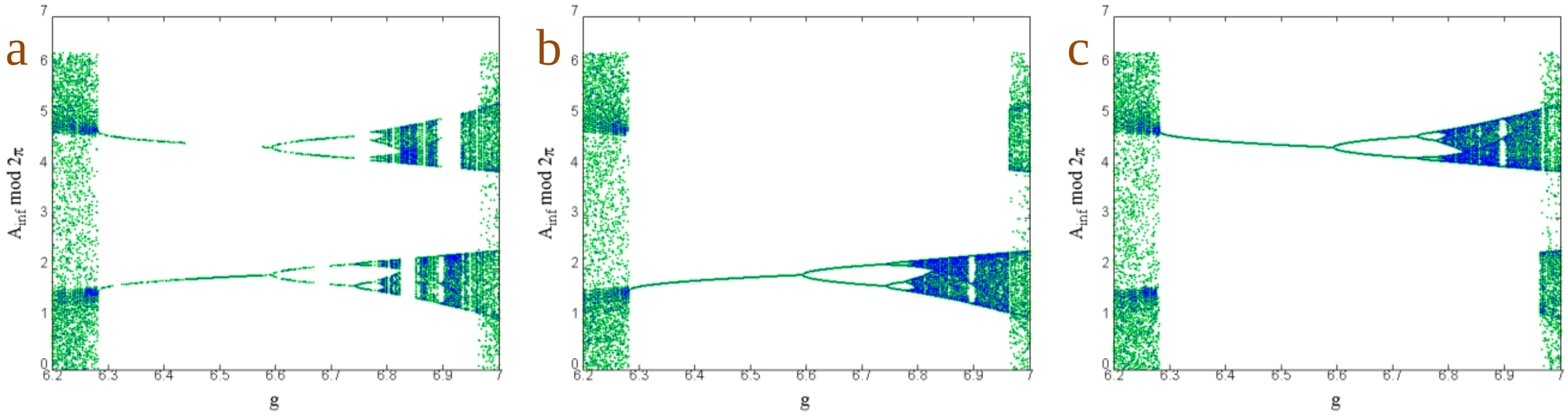

**Fig. 10.** $Mod\ 2\pi$ of the bifurcation diagram of Eq. 1. (a) With random initial condition in the interval $A_0 \in [0,2\pi]$. (b) With initial condition $A_0 = \frac{\pi}{2}$. (c) With initial condition $A_0 = \frac{3\pi}{2}$.

## 2.4. Periodic windows

By increasing the $g$ parameter, another type of behavior can occur which is called bioperiodic. These patterns are very sensitive to initial condition and can model many biological rhythms [1]. Black ellipsoids in Fig. 11 show some of periodic windows in the logarithm of the bifurcation diagram and its $mod\ 2\pi$. Phase portrait and cobweb plot of a bioperiodic pattern in $g = 9.21$ are shown in Fig. 12. As it shows, this periodic behavior crosses through a positive peak and a negative peak in two different periods of sine. It is obvious that these windows occur in the high values of $g$ parameter, thus there is not any stable fixed point in Eq. 1 and these periodic behaviors begin from period two. Fig. 13 represents bifurcation diagram of the first periodic window which is investigated in Fig. 11. This figure shows that period doubling route to chaos happens, then the periodic window has ended and the biotic patterns appear.

In order to investigate the location of the first mentioned periodic window, consider the cobweb plot of Fig. 12. The first parameter which can generate this type of periodic window is the one that has the following conditions,

$$\begin{cases} f^n(A^*) > 3\pi & if\ n\ is\ odd \\ f^n(A^*) < \pi & if\ n\ is\ even \end{cases} . \qquad (6)$$

$A^*$ is the critical point calculated in Eq. 3 and $f^n(.) = \underbrace{f\left(f\left(f \ldots (f(.))\right)\right)}_{n\ times}$.

Part (a) of Fig. 14 investigates the conditions of Eq. 6 until $n = 3$. Three horizontal black lines in the figure show the location of $\pi, 3\pi$ and $5\pi$. Black rectangular shapes depict regions satisfy the conditions of Eq. 6. The first rectangle determines the region $g \in [9.0065, 9.7558]$ for the first mentioned periodic window. In order to find the more precise region of this periodic window, we compare the results of these three functions with the result of $f^{20}$ which is shown in part (b) of

Fig. 14. The comparison investigates $g \in [9.205, 9.6838]$ as the region of this periodic window which is completely matched with the bifurcation diagram of Fig. 13.

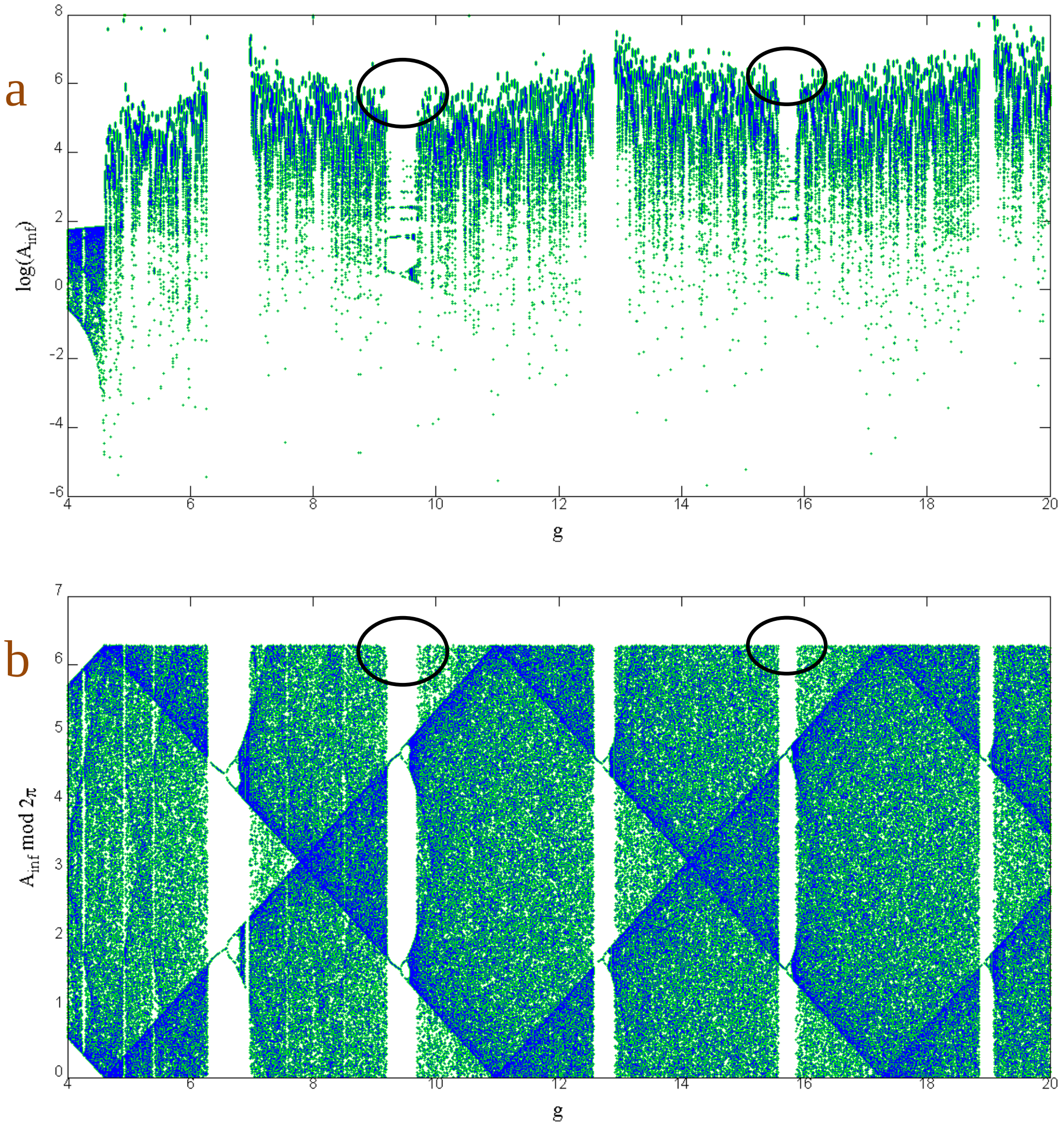

**Fig. 11.** (a) Logarithm of the bifurcation diagram of Eq. 1. Black ellipsoids show periodic windows. (b) $Mod\ 2\pi$ of the bifurcation diagram of Eq. 1. Black ellipsoids show periodic windows.

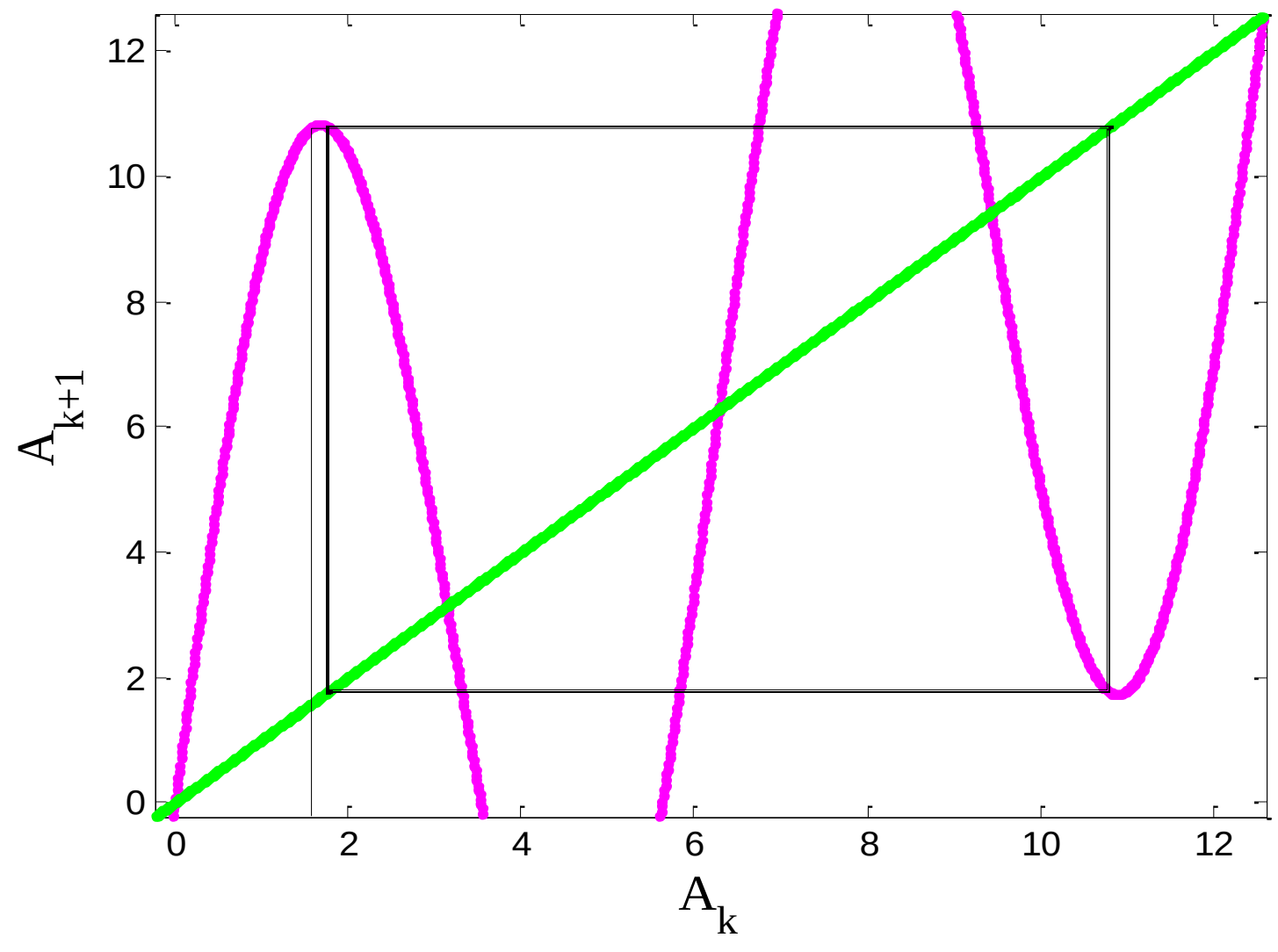

**Fig. 12.** Cobweb plot of a bioperiodic pattern in $g = 9.21$ and initial condition $A_0 = \frac{\pi}{2}$.

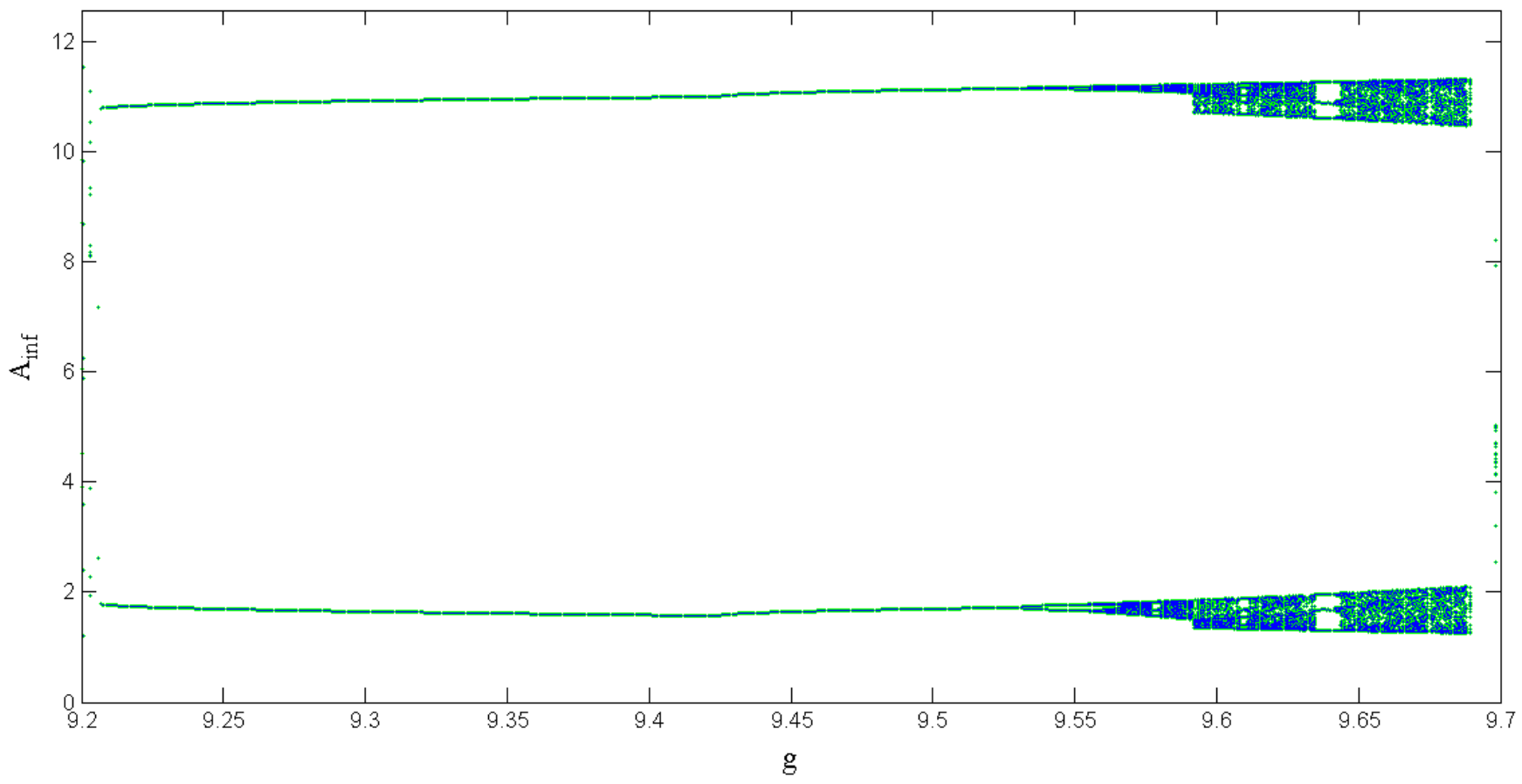

**Fig. 13**. Bifurcation diagram of the first periodic window in Fig. 9 with initial condition $A_0 = \frac{\pi}{2}$.

A later periodic window of this type has the same mechanism in phase space as the previous one with a little difference which is contained three periods of sine as shown in cobweb plot of Fig. 15. It is obvious that the condition which can generate this type of periodic window is as follows,

$$\begin{cases} f^n(A^*) > 5\pi & if\ n\ is\ odd \\ f^n(A^*) < \pi & if\ n\ is\ even \end{cases} \quad (7)$$

Considering part (a) of Fig. 14 shows that the system satisfies these conditions in the range $g \in [15.501, 15.8837]$. As discussed in the previous window, the more precise region of this window which is obtained by a comparison between the results of part (a) and part (c) of Fig. 14 is $g \in [15.5796, 15.8763]$.

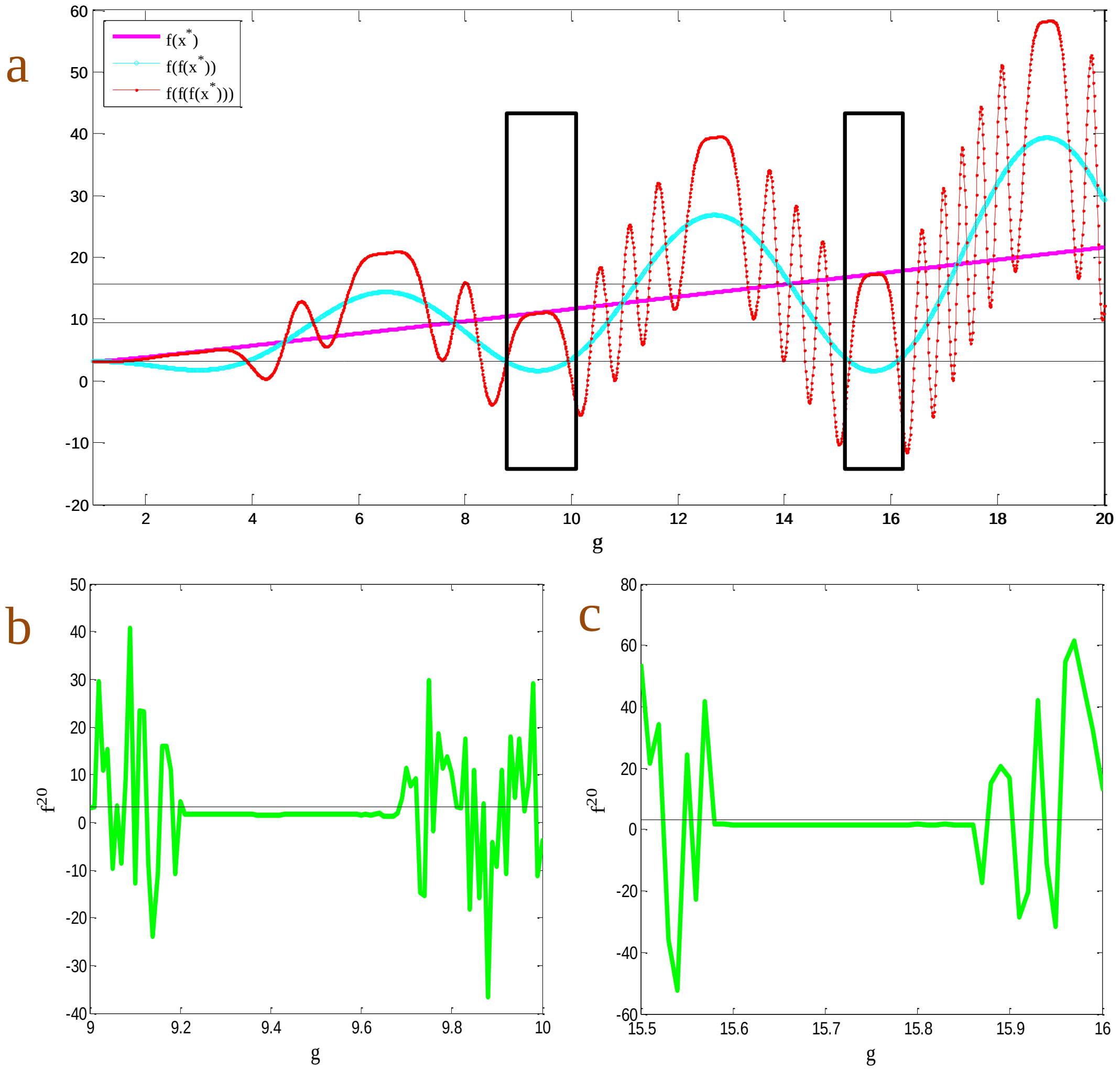

**Fig. 14.** (a) Investigating conditions of Eq. 6 until $n = 3$. Horizontal black lines show the location of $\pi$, $3\pi$ and $5\pi$ and black rectangular shapes depict regions that satisfy conditions of Eq. 6. (b) Investigating conditions of Eq. 6 in $n = 20$ and $g \in [9,10]$. Horizontal black line shows the location of $\pi$. (c) Investigating conditions of Eq. 6 in $n = 20$ and $g \in [15.5,16]$. Horizontal black line shows the location of $\pi$.

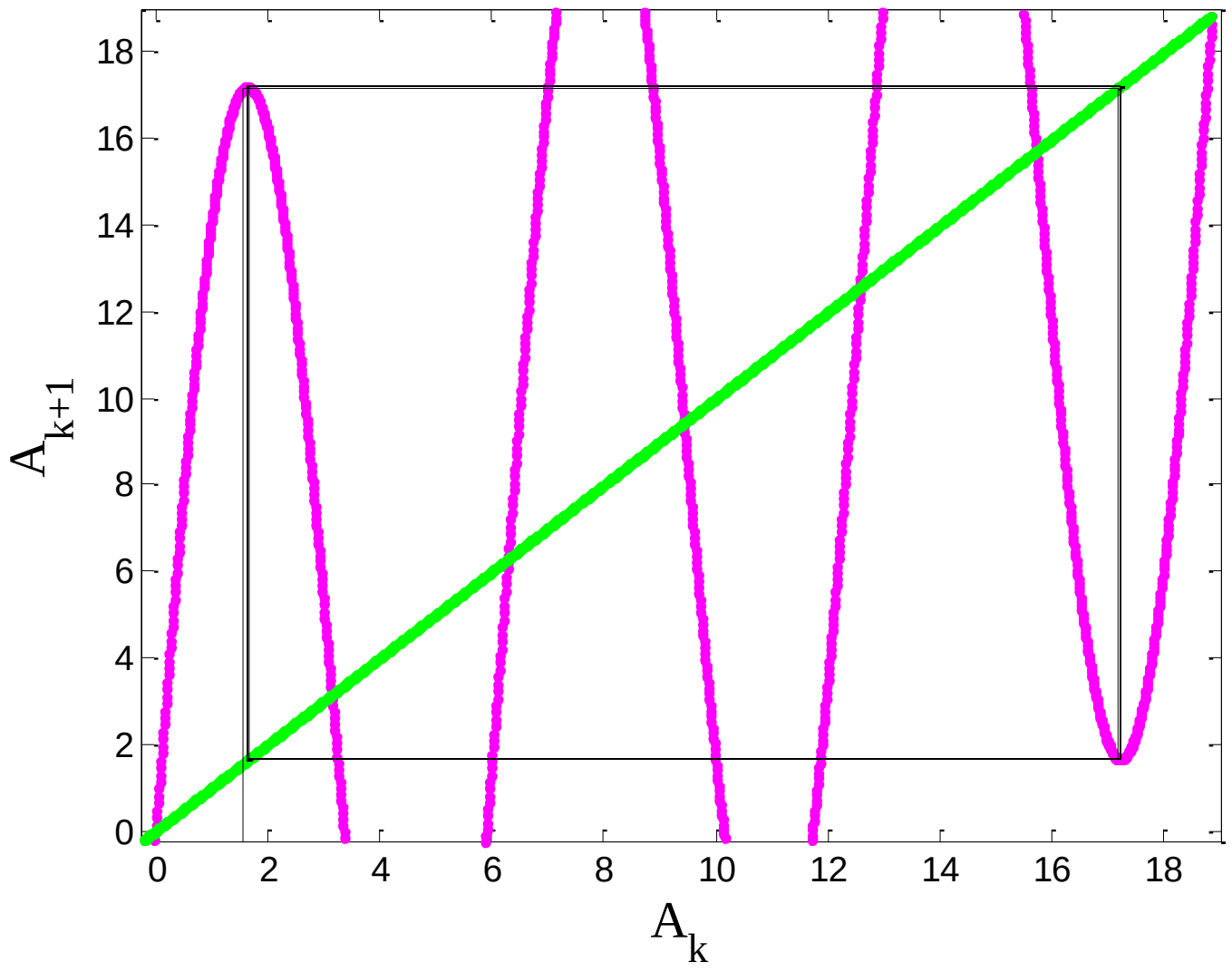

**Fig. 15.** Cobweb plot of a bioperiodic pattern in g $= 15.6$ and initial condition $A_0 = \frac{\pi}{2}$.

It is obvious that the conditions of the second window also satisfy the conditions of the first window. Thus, all of those kinds of periodic windows satisfy conditions of Eq. 6. In other words, the conditions determine the general location of all discussed periodic windows and more details of the location are obtained using the specific conditions for each window.

## 3. Discussion

In this paper, different behaviors of process equation are investigated with a different viewpoint to the phase spaces. Period doubling route to chaos, bios, periodic windows and unstable windows are discussed with a deep look in their phase space and the location of each behavior is determined. A closer look at the $mod\ 2\pi$ of the bifurcation diagram of the process equation shows that there are some unstable small windows and periodic small windows among biotic behaviors which have a different structure from what is discussed in the previous section. Fig. 16 shows some examples of these windows. Red rectangles indicate unstable small windows and black ones indicate periodic small windows. Now a question is raised that what is the difference of these windows with the previous ones which our methods cannot determine the regions of these windows.

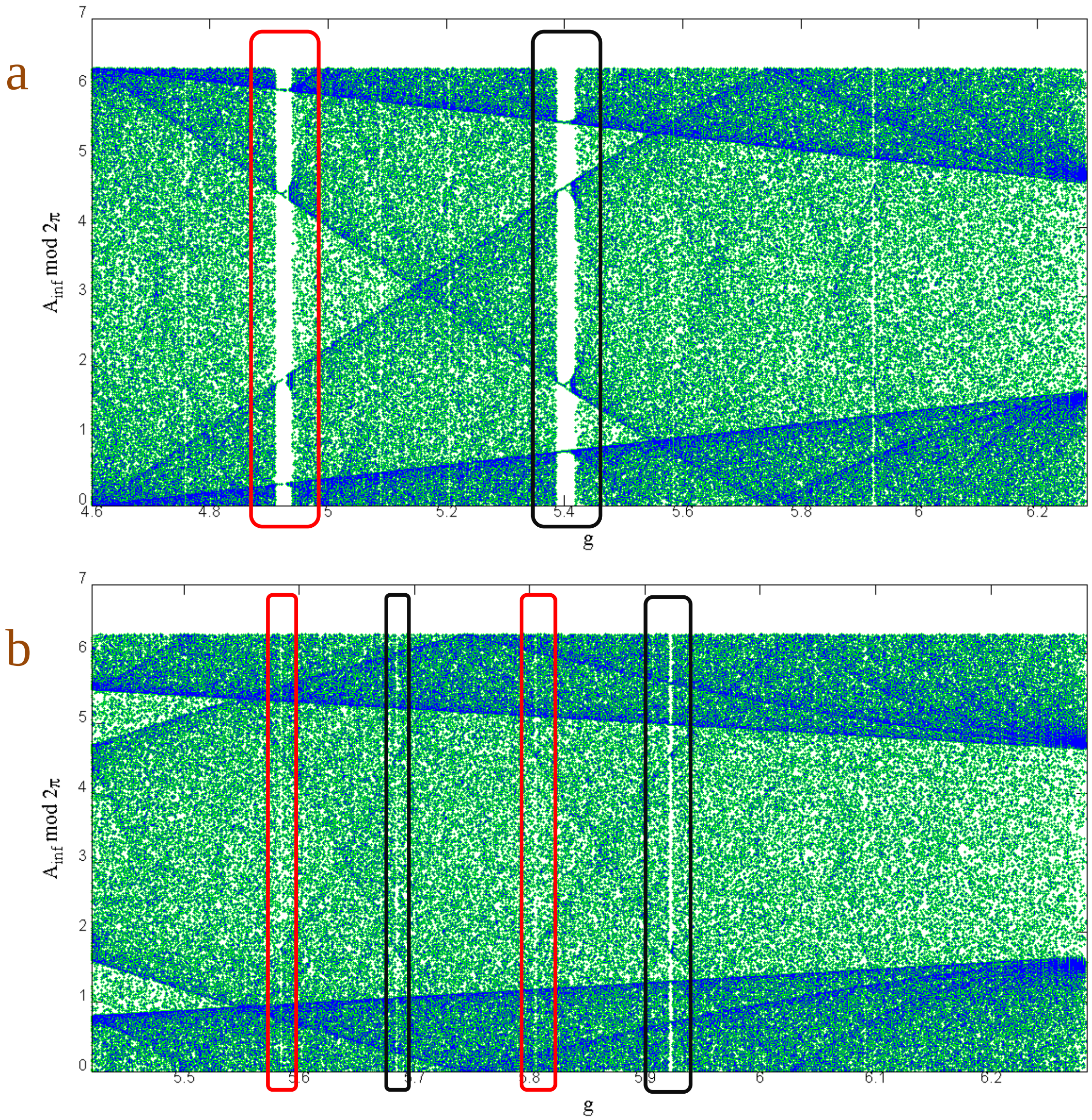

**Fig. 16.** (a) Bifurcation diagram of Eq. 1 in $mod\ 2\pi$ with respect to parameter $g \in [4.6, 2\pi]$ and initial condition $A_0 = \frac{\pi}{2}$. Red rectangles indicate unstable small windows and black ones indicate periodic small windows. (b) Bifurcation diagram of Eq. 1 in $mod\ 2\pi$ with respect to parameter $g \in [5.45, 2\pi]$ and initial condition $A_0 = \frac{\pi}{2}$. Red rectangles indicate unstable small windows and black ones indicate periodic small windows.

In the part (a) of Fig. 16, $mod\ 2\pi$ of the bifurcation diagram is shown with respect to changing $g$ parameter in the interval $g \in [4.6, 2\pi]$. As the figure shows, there are at least two windows in the small region. Part (b) of Fig. 16 magnifies a smaller region, in the interval $g \in [5.45, 2\pi]$, and shows that there are some smaller windows in the small region which just can be seen with a deeper look. In order to detect the difference between these small windows with the previous ones, consider the cobweb plots of Fig. 17. Each part of the figure is an example of the system with a parameter value which is located at the corresponding windows of Fig. 16. In other words, the simplest behavior of each window which is considered in Fig. 16 from left to right and top to down are shown in part (a) to (f) of Fig. 17, respectively.

Parts (a), (c) and (e) of Fig. 17 depict that the structures of these unstable windows are a little different from those that are investigated in section 2.3. From the cobweb plots of Fig. 17, it is obvious that condition of Eq. 5 cannot determine the regions of these windows. Considering the cobweb plots, the beginning of these windows can be calculated using condition of Eq. 8 for different $L$s.

$$A_{k+L} - A_k = 2\pi \tag{8}$$

Assuming $L = 2$, we have

$$A_{k+2} - A_k = 2\pi \rightarrow g\sin(A_k) + g\sin(A_k + g\sin(A_k)) = 2\pi. \tag{9}$$

This condition is more complex than Eq. 5 and we try to solve it with a numerical method. It has many different solutions which depict lots of small unstable windows. One of solutions is shown in the Eq. 10.

$$g = 4.911529662464309878494343949 8103 \quad A_k = 1.7300 \tag{10}$$

This solution has the same behavior as is shown in part (a) of Fig. 17. Bifurcation diagram of Fig. 18 depicts that the two behaviors are parts of a same unstable window. Another attractive subject in these unstable small windows is the occurrence of odd periods at the beginning of these windows as shown in part (c) and (e) of Fig. 17 (period three). Those behaviors satisfy condition of Eq. 8 with $L = 3$.

Comparing cobweb plots of parts (b), (d) and (f) of Fig. 17 with cobweb plots of previous periodic windows in Fig. 12 and Fig. 15 clarifies the structural difference of these two kinds of windows. In the previous periodic windows, the state of the system jumps from one period of sine to another one in an opposite sign of sine phase just in one iteration, but in the new windows the periodic behavior starts from period four or six and the state jumps from one period of sine to another period with the same sign and goes to the opposite sign with two or more steps. The parameter which can generate this type of periodic window is the one that has the following conditions,

$$\begin{cases} f^{Ln}(A^*) > (2I+1)\pi & if\ n\ is\ odd \\ f^{Ln}(A^*) < \pi & if\ n\ is\ even \end{cases}. \tag{11}$$

This condition can be investigated with different $L$ and $I$ parameters. Part (b) of Fig. 17 satisfies condition of Eq. 11 with $L = 2$ and $I = 1$, part (d) satisfies this condition with $L = 3$ and $I = 1$, and finally part (f) satisfies this condition with $L = 3$ and $I = 2$. Fig. 19 investigates the condition of Eq. 11 for windows containing the behaviors of part (b), (d) and (f) of Fig. 17. As the figure shows, there are lots of other windows that satisfies the condition of Eq. 11 and thus there are many small windows among the biotic behavior.

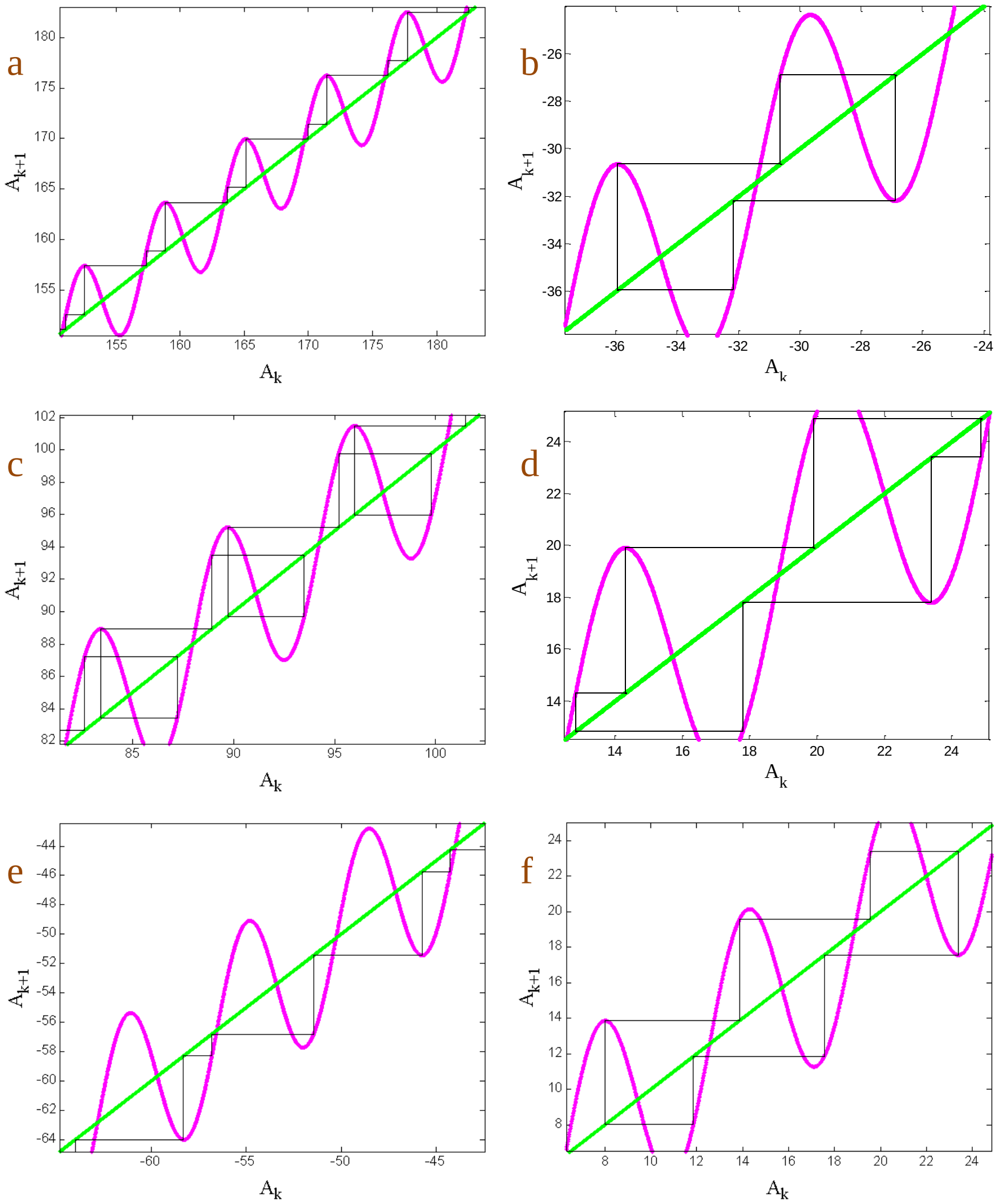

**Fig. 17.** Cobweb plot of termination of signal after passing the transient time. (a) With parameter $g = 4.915$ and initial condition $A_0 = \frac{\pi}{2}$. (b) With parameter $g = 5.39$ and initial condition $A_0 = \frac{\pi}{2}$. (c) With parameter $g = 5.5835$ and initial condition $A_0 = \frac{\pi}{2}$. (d) With parameter $g = 5.6845$ and initial condition $A_0 = \frac{\pi}{2}$. (e) With parameter $g = 5.805$ and initial condition $A_0 = \frac{\pi}{2}$. (f) With parameter $g = 5.921$ and initial condition $A_0 = \frac{\pi}{2}$.

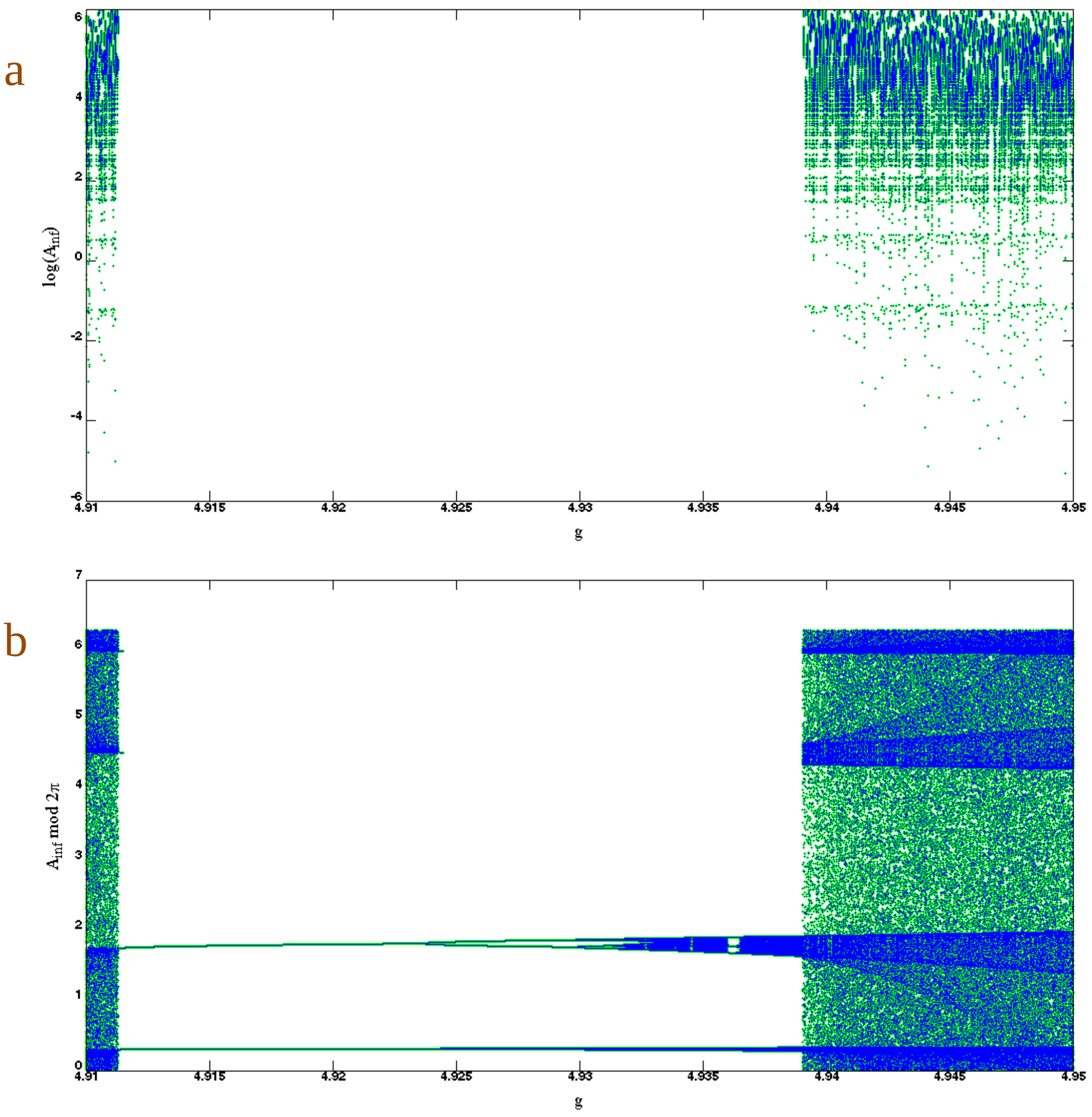

**Fig. 18.** Bifurcation diagram of Eq. 1 with respect to parameter $g \in [4.91,4.95]$ and initial condition $A_0 = 1.7300$. (a) Logarithm of bifurcation diagram (b) $mod\ 2\pi$ of bifurcation diagram.

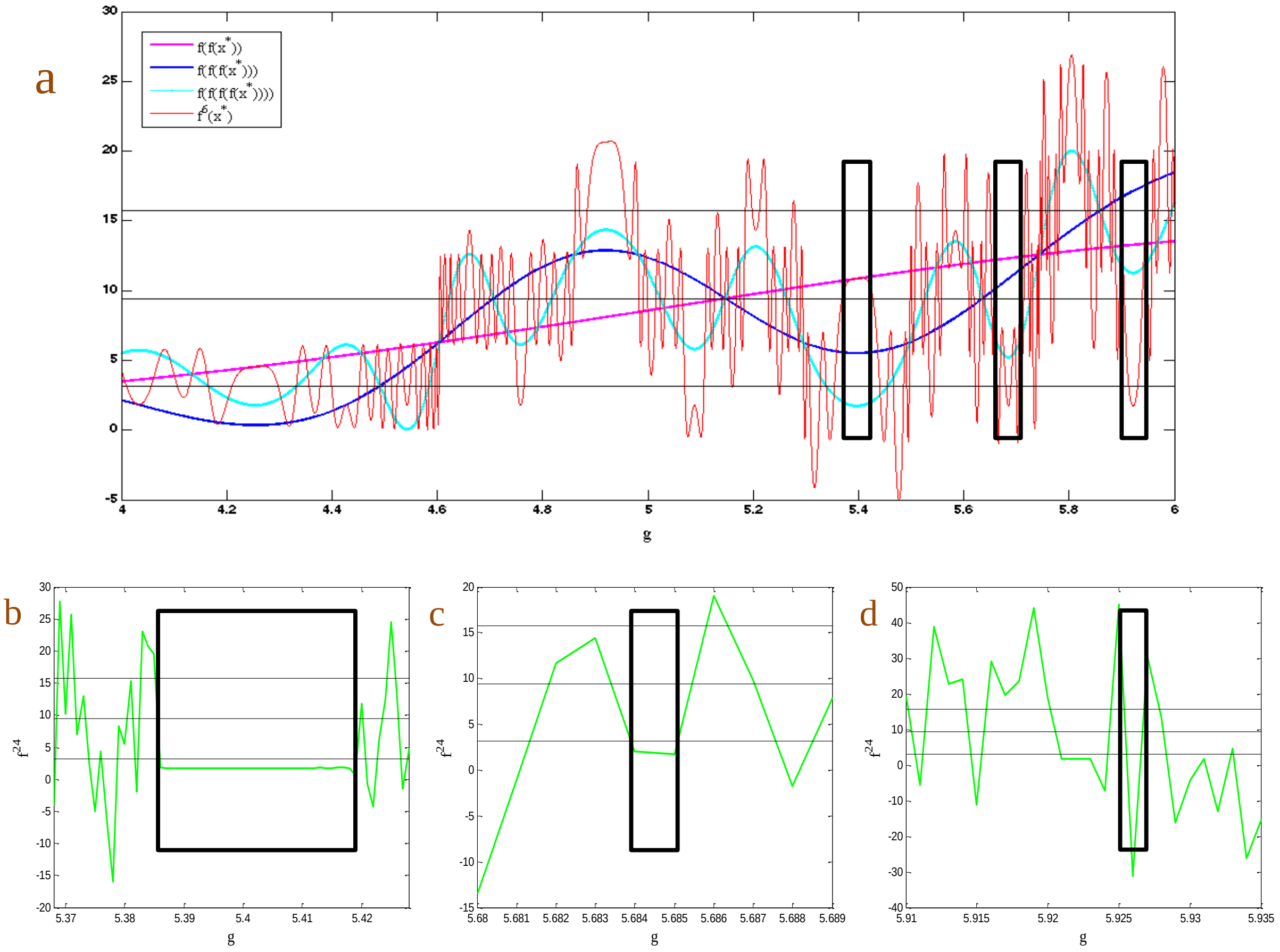

**Fig. 19.** (a) Investigating conditions of Eq. 11 until n = 6. Horizontal black lines show the location of π, 3π and 5π and black rectangular shapes depict regions that satisfy conditions of Eq. 11 in the range of parameters of the cobweb plots of part (b), (d) and (f) of Fig. 17. (b) Investigating conditions of Eq. 11 for part (b) of Fig. 17 with $L = 2, I = 1$, n = 24 and $g \in [5.3683, 5.44]$. (c) Investigating conditions of Eq. 11 for part (d) of Fig. 17 with $L = 3, I = 1$, n = 24 and $g \in [5.68, 5.689]$. (d) Investigating conditions of Eq. 11 for part (f) of Fig. 17 with $L = 3, I = 2$, n = 24 and $g \in [5.91, 5.935]$.

Fig. 20 shows Q-curves of process equations. Q curves are polynomial curves which show the dense parts of bifurcation diagram [28]. The process equation in each period of sine has two critical points as follows,

$$f'(A) = 1 + g\cos(A) \rightarrow f' = 0 \rightarrow A_1^* = \cos^{-1}\left(-\frac{1}{g}\right) \;\&\; A_2^* = 2\pi - \cos^{-1}\left(-\frac{1}{g}\right). \quad (12)$$

Four Q-curves, $Q_1(A_1^*), Q_2(A_1^*), Q_1(A_2^*), Q_2(A_2^*)$, are shown in Fig. 20 with different colors. These Q-curves express the hypothesis of a difference between the windows discussed in section 2 and the new windows because of the difference in their dense lines of the bifurcation diagram in $mod\ 2\pi$.

.

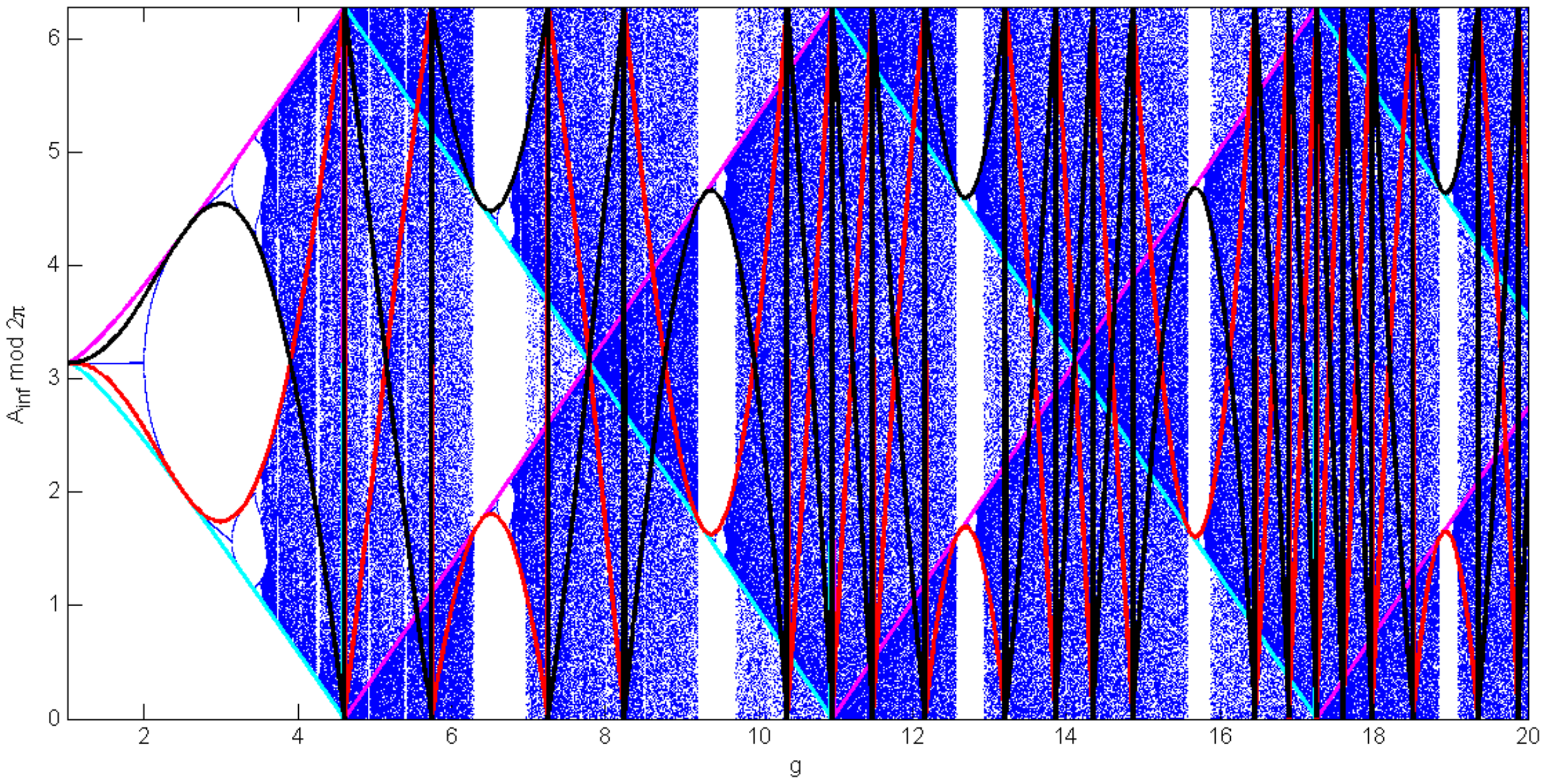

**Fig. 20.** Bifurcation diagram of Eq.1 in $mod\ 2\pi$ and its Four Q-curves.

It is obvious that all calculations of section 2 are based on $A_1^*$ and the the results of $A_2^*$ are absolutely the same as $A_1^*$.

## 5. Conclusion

This paper focuses on the different behaviors of process equation as a multistable one dimensional map with nonlinear feedback, and the parameters of their occurrences are calculated. Using a deeper look to phase portrait, cobweb plot and higher iterate of process equation, period doubling route to chaos of the equation are investigated and unifurcations are discussed. Then, the parameter of process equation’s entrance to biotic pattern and some features of the behavior are explained with looking to its phase portrait. In this paper, various patterns of unstable and periodic windows are shown and the reasons of their happenings are investigated in details. Also, we discuss about some other types of unstable and periodic windows and the reason of their happenings. On the other hands, Q-curves determine the difference of these windows with the previous ones.

## Acknowledgements

We dedicate this paper to the memory of Hector Sabelli.